\documentclass[11pt, reqno]{amsart}
 
 \usepackage{amsmath}
 \usepackage{amssymb}
\usepackage{bm}

\newcommand{\mysection}[1]{\section{#1}
      \setcounter{equation}{0}}

\newtheorem{theorem}{Theorem}[section]
\newtheorem{lemma}[theorem]{Lemma}

\newtheorem{corollary}[theorem]{Corollary} 

\theoremstyle{definition}
\newtheorem{assumption}{Assumption}[section]

\theoremstyle{remark}
\newtheorem{remark}{Remark}[section]

\newcommand{\tr}{\text{\rm tr}\,}

\newcommand\bbeta{\text{\raise-.2ex\hbox{$\bm{\beta}$}}}

\newcommand\bR{\mathbb{R}}

\newcommand\bB{\mathbb{B}}

\newcommand\cF{\mathcal{F}}

\newcommand\frA{\mathfrak{A}}
\newcommand\frB{\mathfrak{B}}

\newcommand\infsup{\operatornamewithlimits{inf\,\,\,sup}}
\newcommand\supinf{\operatornamewithlimits{sup\,\,\,inf}}
\newcommand\supsup{\operatornamewithlimits{sup\,\,\,sup}}

\begin{document}

\title[Approximating the value functions]
{Approximating the value functions
for stochastic differential games   with the ones
having bounded second derivatives}

\author{N.V. Krylov}
\thanks{The  author was partially supported by
 NSF Grant DMS-1160569}
\email{krylov@math.umn.edu}
\address{127 Vincent Hall, University of Minnesota,
 Minneapolis, MN, 55455}



\begin{abstract}
We show a method of uniform approximation of
 the value functions of
uniformly nondegenerate stochastic differential games
in smooth domains   up to
a constant over $K$
with the ones having second-order derivatives
bounded by a constant times $K$ for any $K\geq1$.
 
\end{abstract}

\maketitle

\mysection{Introduction}

Let $\bR^{d}=\{x=(x^{1},...,x^{d})\}$
be a $d$-dimensional Euclidean space and  let $d_{1}\geq d$
  be an integer.
Assume that we are given separable metric spaces
  $A$ and $B$,   and let,
for each $\alpha\in A$, $\beta\in B$, 
  the following 
  functions on $\bR^{d}$ be given: 

(i) $d\times d_{1}$
matrix-valued $\sigma^{\alpha\beta}( x)
 =
(\sigma^{\alpha\beta}_{ij}( x))$,

(ii)
$\bR^{d}$-valued $b^{\alpha\beta}( x)=
(b^{\alpha\beta}_{i }( x))$, and

(iii)
real-valued  functions 
$c^{\alpha\beta}( x)\geq0 $,   
  $f^{\alpha\beta}( x) $, and  
$g(x)$. 

Under natural assumptions which will be specified later
one   associates with these objects and a bounded domain 
$G\subset\bR^{d}$
a stochastic differential
game with the diffusion term $\sigma^{\alpha\beta}(x)$,
  drift term $b^{\alpha\beta}(x)$, discount rate 
$c^{\alpha\beta}(x)$, running cost $f^{\alpha\beta}(x)$,
and the final cost $g(x)$ payed when the underlying process
first exits from $G$.

After the order of players is specified in a certain way 
it turns 
out (see, for instance, \cite{FS89},
 \cite{Ko09}, \cite{Sw96} or Remark 2.2 in \cite{Kr_14}) 
that the value function $v(x)$
of this differential game is a unique 
continuous in $\bar{G}$ viscosity
solution
of the Isaacs equation
$$
H[v]=0
$$
in $G$ with boundary condition $v=g$ on $\partial G$,
where for a sufficiently smooth function $u=u(x)$  
$$
L^{\alpha\beta} u(x):=a^{\alpha\beta}_{ij}( x)D_{ij}u(x)+
b ^{\alpha\beta}_{i }( x)D_{i}u(x)-c^{\alpha\beta} ( x)u(x),
$$
$$
a^{\alpha\beta}( x):=(1/2)\sigma^{\alpha\beta}( x)
(\sigma^{\alpha\beta}( x))^{*},\quad
D_{i}=\partial/\partial x^{i},\quad D_{ij}=D_{i}D_{j},
 $$
\begin{equation}
                                                     \label{1.16.1}
H[u](x)=\supinf_{\alpha\in A\,\,\beta\in B}
[L^{\alpha\beta} u(x)+f^{\alpha\beta} (x)].
\end{equation}

Under some assumptions
one explicitly constructs a convex positive-ho\-mo\-ge\-ne\-ous
of degree one function $P(u_{ij},u_{i},u)$ such that for any $K\geq1$
the equation
\begin{equation}
                                               \label{3.13.1}
\max(H[u],P[u]-K)=0
\end{equation}
in $G$ with boundary condition $u=g$ on $\partial G$
has a unique solution $v_{K}$ in class $C^{1,1}_{loc}(G)
\cap C(\bar{G})$ with the second-order
derivatives bounded by a constant times $K$ divided
by the distance to the boundary. Here
$$
P[u](x)=P(D_{ij}u(x),D_{i}u(x),u(x)).
$$

The goal of this article is to prove
the conjecture stated in \cite{Kr_12}:
$|v-v_{K}|\leq N/K$ in $G$ for $K\geq1$, where $N$
is independent of $K$. Such a result even in a much weaker form
was already used in numerical approximations
of solutions of the Isaacs equations in \cite{Kr_14_2}.

The result belongs to the theory of partial differential
 equations.
However, the proof  we give is purely
probabilistic and quite nontrivial involving,
in particular, a reduction of differential
games in domains to the ones on a smooth manifolds
without boundary. The main idea  underlying
this reduction is explained in the last two
sections of \cite{Kr95}
and, of course, we represent $v_{K}$
also as a value function for a 
corresponding stochastic differential game.
Still it is worth mentioning that the methods
of the theory of partial differential
 equations can be used to obtain results similar
to ours albeit not that sharp in what concerns the rate
of approximations even though
for Isaacs equations with much less regular coefficients
than ours
(see \cite{Kr_15}).

The article is organized as follows.
Section \ref{section 2.26.3} contains our main result.
In Section \ref{section 4.9.1} we prove the dynamic
programming principle for stochastic differential games in the whole
space. In Section \ref{section 2.27.1} we show how to
reduce the stochastic differential game
in a domain to the one in the whole space
having four more dimensions. Actually,
the resulting stochastic differential games lives
on a closed manifold without boundary.
In the final Section \ref{section 4.9.2} we prove
our main result, Theorem \ref{theorem 1.14.01}.
 
\mysection{Main result }
                                    \label{section 2.26.3}

We start with our assumptions.
 
\begin{assumption}
                                    \label{assumption 1.9.1}

(i) The functions $\sigma^{\alpha\beta}(x)$, $b^{\alpha\beta}(x)$,
$c^{\alpha\beta}(x)$, and $f^{\alpha\beta}(x)$
are continuous with respect to
$\beta\in B$ for each $(\alpha, x)$ and continuous with respect
to $\alpha\in A$ uniformly with respect to $\beta\in B$
for each $x$. 

(ii) 
for any $x \in\bR^{d}$ 
 and $(\alpha,\beta )\in A\times B $
$$
  \|\sigma^{\alpha\beta}( x )\|,|b^{\alpha\beta}( x )|,
|c^{\alpha\beta}( x )|,|f^{\alpha\beta}( x )|
\leq K_{0},
$$
where   $K_{0}$ is
a fixed constants and for a matrix $\sigma$ we denote $\|\sigma\|^{2}
=\tr\sigma\sigma^{*}$,

(iii) For any $(\alpha,\beta )\in A\times B $ and $x,y\in\bR^{d}$
we have
$$
\|\sigma^{\alpha\beta}( x )-\sigma^{\alpha\beta}( y )\|,
|u^{\alpha\beta}( x )-u^{\alpha\beta}( y )| \leq K_{0}|x-y|,
$$
where $u=b,c,f$.
\end{assumption}

Let $(\Omega,\cF,P)$ be a complete probability space,
let $\{\cF_{t},t\geq0\}$ be an increasing filtration  
of $\sigma$-fields $\cF_{t}\subset \cF $ such that
each $\cF_{t}$ is complete with respect to $\cF,P$.
We suppose that on $(\Omega,\cF,P)$ we are given a
$d_{1}$-dimensional
  Wiener processes $w_{t}$,
which is a Wiener processes relative to $\{\cF_{t}\}$.

 The set of progressively measurable $A$-valued
processes $\alpha_{t}=\alpha_{t}(\omega)$ is denoted by $\frA$. 
Similarly we define $\frB$
as the set of $B$-valued  progressively measurable functions.
By  $ \bB $ we denote
the set of $\frB$-valued functions 
$ \bbeta(\alpha_{\cdot})$ on $\frA$
such that, for any $T\in(0,\infty)$ and any $\alpha^{1}_{\cdot},
\alpha^{2}_{\cdot}\in\frA$ satisfying
\begin{equation}
                                                  \label{4.5.4}
P(  \alpha^{1}_{t}=\alpha^{2}_{t} 
 \quad\text{for almost all}\quad t\leq T)=1,
\end{equation}
we have
$$
P(  \bbeta_{t}(\alpha^{1}_{\cdot})=\bbeta_{t}(\alpha^{2}_{\cdot}) 
\quad\text{for almost all}\quad t\leq T)=1.
$$

Fix a domain $G\subset\bR^{d}$, and impose the following.

\begin{assumption}
                                 \label{assumption 3.19.1}
 $G$ is a bounded domain of class $C^{3}$,
$g\in C^{3}$,
and there exists a constant $\delta\in(0,1)$ such that
for any $\alpha\in A$, $\beta\in B$,  
and $x,\lambda\in \bR^{d}$
$$
\delta|\lambda|^{2}\leq a^{\alpha\beta}_{ij}( x)
\lambda^{i}\lambda^{j}\leq \delta^{-1}|\lambda|^{2}.
$$

\end{assumption}

\begin{remark}
                                                     \label{remark 2.22.2}

As is well known, if Assumption \ref{assumption 3.19.1} 
is satisfied, then
 there exists a bounded from above
  $\Psi\in   C^{3}_{loc}(\bR^{d})$
such that 
  $\Psi>0$ in $G$, $\Psi=  0$ and $|D\Psi|\geq1$
 on $\partial G$, $\Psi(x)\to-\infty$ as
$|x|\to\infty$, and
for all $\alpha\in A$, $\beta\in B$,  
and $x\in G$
\begin{equation}
                                       \label{3.20.1}
  L^{\alpha\beta}\Psi( x)
+c^{\alpha\beta}\Psi( x)\leq-1.
\end{equation}
\end{remark}

For $\alpha_{\cdot}\in\frA$, 
$ \beta_{\cdot} 
\in\frB$, and $x\in\bR^{d}$ consider the following It\^o
equation
\begin{equation}
                                             \label{5.11.1}
x_{t}=x+\int_{0}^{t}
\sigma^{\alpha_{s}\beta_{s} }(  x_{s})
\,dw _{s}
+\int_{0}^{t} b^{\alpha_{s}
\beta_{s} }(  x_{s}) \,ds.
\end{equation}
Observe that for any $\alpha_{\cdot}\in\frA$, 
$ \beta_{\cdot} 
\in\frB$,  $x\in\bR^{d}$, and $T\in(0,\infty)$ it has a unique solution
on $[0,T]$
which we denote by $x^{\alpha_{\cdot} 
\beta_{\cdot} x}_{t} $.

Set
$$
\phi^{\alpha_{\cdot}\beta_{\cdot} x}_{t}
=\int_{0}^{t} 
c^{\alpha_{s}
\beta_{s} }( x^{\alpha_{\cdot}\beta_{\cdot}  x}_{s}) 
\,ds,
$$
 define $\tau^{\alpha_{\cdot}\beta_{\cdot} x}$ as the first exit
time of $x^{\alpha_{\cdot} 
\beta_{\cdot} x}_{t}$ from $G$, and introduce
\begin{equation}
                                                    \label{2.12.2}
v(x)=\infsup_{\bbeta\in\bB\,\,\alpha_{\cdot}\in\frA}
E_{x}^{\alpha_{\cdot}\bbeta(\alpha_{\cdot})}
\big[\int_{0}^{\tau}
 f ( x_{t})e^{-\phi_{t}
 }\,dt+g(x_{\tau})e^{-\phi_{\tau}
 }\big],
\end{equation}
where, as usual,
 the indices $\alpha_{\cdot}$, $\bbeta$, and $x$
at the expectation sign are written  to mean that
they should be placed inside the expectation sign
wherever and as appropriate, for instance,
$$
E_{x}^{\alpha_{\cdot}\beta_{\cdot}}\big[\int_{0}^{\tau}
f( x_{t})e^{-\phi_{t} }\,dt+g(x_{\tau})e^{-\phi_{\tau}
 }\big]
$$
$$
:=
E \big[\int_{0}^{\tau^{\alpha_{\cdot}\beta_{\cdot}  x}}
f^{\alpha_{t}\beta_{t}   }
( x^{\alpha_{\cdot}\beta_{\cdot}  x}_{t})
e^{-\phi^{\alpha_{\cdot}\beta_{\cdot}  x}_{t}
 }\,dt
+g(x^{\alpha_{\cdot}\beta_{\cdot}  x}
_{\tau^{\alpha_{\cdot}\beta_{\cdot}  x}}
)
e^{-\phi^{\alpha_{\cdot}\beta_{\cdot}  x}
_{\tau^{\alpha_{\cdot}\beta_{\cdot}  x}}}\big].
$$

Observe that, formally, the value $x_{\tau}$ 
 may not be defined if $\tau=\infty$.
In that case we set the corresponding terms to equal zero.
 This is natural because It\^o's formula easily yields
that $E^{\alpha_{\cdot}\beta_{\cdot}}_{x}\tau\leq \Psi(x)$
in $G$, so that $\tau<\infty$ (a.s.).

We also need a few new objects.  
In the end of Section 1 of \cite{Kr_12} a 
function $P(u_{ij},u_{i},u)$ is constructed
defined for all symmetric $d\times d$ matrices $(u_{ij})$,
$\bR^{d}$-vectors $(u_{i})$, and $u\in\bR$ such that
  it is positive-homogeneous of degree
one, is
Lipschitz continuous, and at all points of differentiability of $P$
for all values of arguments we have $P_{u}\leq0$ and
$$
\hat\delta|\lambda|^{2}\leq P_{u_{ij}}
\lambda^{i}\lambda^{j}\leq \hat{\delta}^{-1}|\lambda|^{2},
$$
where $\hat{\delta}$ is a constant in $(0,1)$
depending only on
$d,K_{0}$, and $\delta$.  

We now state a part of Theorem 1.1 of \cite{Kr_12}
which we need.
\begin{theorem}
                                    \label{theorem 9.23.1}
For any $K\geq0$   the equation
\begin{equation}
                                               \label{9.23.2}
\max(H[u],P[u]-K)=0
\end{equation}
in $G$  (a.e.) with boundary condition $v=g$ on $\partial G$
has a unique solution 
$u\in C^{0,1}(\bar{G})\cap C^{1,1}_{loc}(G)$.
\end{theorem}

Our  main result consists of proving
the conjecture stated in \cite{Kr_12}.

\begin{theorem}
                                         \label{theorem 1.14.01}
 
Denote by $u_{K}$ the function from 
Theorem \ref{theorem 9.23.1}.
Then there exists a constant $N$
such that
$|v-u_{K}|\leq N\Psi/K$ in $G$ for $K\geq1$.
\end{theorem}

\mysection{On degenerate stochastic differential games
in the whole space}
                                          \label{section 4.9.1}

Here we suppose that the assumptions of Section
\ref{section 2.26.3} are satisfied with the following exceptions.
We
do not need Assumption \ref{assumption 1.9.1} (iii)
satisfied for $u=c,f$. It suffices to have
the functions 
$c^{\alpha\beta}( x )$
and
$f^{\alpha\beta}( x )$
  uniformly continuous with
respect to $x$ uniformly with respect to 
$(\alpha,\beta )\in A\times B $. We also abandon
Assumption \ref{assumption 3.19.1}
regarding $G$ and the uniform nondegeneracy
of $a$, but impose the following.

\begin{assumption}
                               \label{assumption 2.23.03}
There exists a constant $\delta_{1}>0$ such that
for any $\alpha\in A$, $\beta\in B$,  
and $x \in \bR^{d}$  
$$
 c^{\alpha\beta} ( x)\geq \delta_{1}.
$$
\end{assumption}

The probability space here and the underlying filtration
of $\sigma$-fields are not necessarily the same
as in Section \ref{section 2.26.3} and in our applications
they indeed will be different. Therefore,
the following assumption is harmless for the purpose
of our applications.
\begin{assumption}
                               \label{assumption 2.23.3}
There exists a $d$-dimensional Wiener
process $\bar{w}_{\cdot}$ which is
a Wiener process relative to
$\{\cF_{t}\}$ and is independent of $w_{\cdot}$.
\end{assumption}

We also use a somewhat different definition
of $v(x)$. Set
$$
v(x)=\infsup_{\bbeta\in\bB\,\,\alpha_{\cdot}\in\frA}
E_{x}^{\alpha_{\cdot}\bbeta(\alpha_{\cdot})}
 \int_{0}^{\infty}
 f ( x_{t})e^{-\phi_{t}
 }\,dt.
$$

The goal of this section is to present the 
proof of the following dynamic programming principle.

\begin{theorem}
                                            \label{theorem 1.14.1}
 
Under the above assumptions

(i) The function $v(x)$   is bounded and 
uniformly continuous in $\bR^{d}$.

(ii)
Let   $\gamma^{\alpha_{\cdot}\beta_{\cdot}x} $
 be an $\{\cF_{t}\}$-stopping
time defined for each $\alpha_{\cdot}\in\frA$, $\beta_{\cdot}\in\frB$,
and $x\in \bR^{d}$. Also let $\lambda_{t}^{\alpha_{\cdot}
\beta_{\cdot}x}\geq0$ be progressively measurable  functions on $\Omega
\times[0,\infty)$ 
defined for each $\alpha_{\cdot}\in\frA$, $\beta_{\cdot}\in\frB$,
and $x\in \bR^{d}$ and
such that they have finite integrals over finite time intervals
(for any $\omega$).
Then for any $x$
\begin{equation}
                                                  \label{1.14.1}
v(x)=\infsup_{\bbeta\in\bB\,\,\alpha_{\cdot}\in\frA}
E_{x}^{\alpha_{\cdot}\bbeta(\alpha_{\cdot})}\big[
v(x_{\gamma})e^{-\phi_{\gamma}-\psi_{\gamma}}
+\int_{0}^{\gamma}
\{f( x_{t})+\lambda_{t}v(x_{t})\}e^{-\phi_{t}-\psi_{t}}\,dt \big],
\end{equation}
where inside the expectation sign
$\gamma=\gamma^{\alpha_{\cdot}\bbeta(\alpha_{\cdot})x}
$ and
$$
\psi^{\alpha_{\cdot}\beta_{\cdot} x}_{t}
=\int_{0}^{t}
\lambda^{\alpha_{\cdot}\beta_{\cdot} x}_{s}\,ds.
$$ 
\end{theorem}

Proof.
For $\varepsilon>0$, $\alpha_{\cdot}\in\frA$, 
$ \beta_{\cdot} 
\in\frB$, and $x\in\bR^{d}$ denote by 
$x^{\alpha_{\cdot}\beta_{\cdot}x}_{t}(\varepsilon)$
the solution of the equation
$$
x_{t}=x+\varepsilon \bar w _{t}+\int_{0}^{t}
\sigma^{\alpha_{s}\beta_{s} }(  x_{s})
\,dw _{s}
+\int_{0}^{t} b^{\alpha_{s}
\beta_{s} }(  x_{s}) \,ds.
$$
Since the coefficients of these equations
satisfy the global Lipschitz condition,
well-known results about It\^o's equations
imply that there is a constant $N$, depending
only on $K_{0}$, such that
for any $\varepsilon>0$, $\alpha_{\cdot}\in\frA$, 
$ \beta_{\cdot} 
\in\frB$, $T\in(0,\infty)$, and $x\in\bR^{d}$
$$
E_{x}^{\alpha_{\cdot} \beta_{\cdot}}
\sup_{t\leq T}|x_{t}-x_{t}(\varepsilon)|^{2}
\leq N\varepsilon^{2} e^{NT}.
$$
It follows that for any $T\in(0,\infty)$ and $\kappa>0$
\begin{equation}
                                               \label{2.24.1}
\lim_{\varepsilon\downarrow0}\sup_{x\in\bR^{d}}
\supsup_{\alpha_{\cdot}\in\frA\,\,\,\beta_{\cdot} 
\in\frB}
P_{x}^{\alpha_{\cdot} \beta_{\cdot}}
(\sup_{t\leq T}|x_{t}-x_{t}(\varepsilon)|\geq\kappa)=0,
\end{equation}
where the indices $\alpha_{\cdot}, \beta_{\cdot}$, and $x$
at the probability sign act in the same way as at the
expectation sign.

Set
$$
v^{\varepsilon}(x) 
=\infsup_{\bbeta\in\bB\,\,\alpha_{\cdot}\in\frA}
E_{x}^{\alpha_{\cdot}\bbeta(\alpha_{\cdot})}
 \int_{0}^{\infty}
 f^{\alpha_{t}\beta_{t}}( x_{t}(\varepsilon))e^{-\phi_{t}(\varepsilon)
 }\,dt,
$$
where
$$
\phi^{\alpha_{\cdot}\beta_{\cdot} x}_{t}(\varepsilon)
=\int_{0}^{t} 
c^{\alpha_{s}
\beta_{s} }( x^{\alpha_{\cdot}\beta_{\cdot}  x}_{s}(\varepsilon)) 
\,ds.
$$
Observe that
$$
|v(x)-v^{\varepsilon}(x)|\leq
\supsup_{\alpha_{\cdot}\in\frA\,\,\,\beta_{\cdot} 
\in\frB}E_{x}^{\alpha_{\cdot} \beta_{\cdot}}
\int_{0}^{\infty}\big[|
f^{\alpha_{t}\beta_{t}}( x_{t}(\varepsilon))
-f^{\alpha_{t}\beta_{t}}( x_{t} )|e^{-\delta_{1}t}
$$
\begin{equation}
                                       \label{2.24.5}
+K_{0}e^{-\delta_{1}t}\int_{0}^{t}
|
c^{\alpha_{s}\beta_{s}}( x_{s}(\varepsilon))
-c^{\alpha_{s}\beta_{s}}( x_{s} )|\,ds\big]\,dt,
\end{equation}
which owing to \eqref{2.24.1}
and the uniform continuity of $c^{\alpha\beta}(x)$
and $f^{\alpha\beta}(x)$ with respect to $x$ implies that
\begin{equation}
                                               \label{2.24.2}
\lim_{\varepsilon\downarrow0}\sup_{\bR^{d}}
|v^{\varepsilon}-v|=0.
\end{equation}
Next, it is also well known that
there is a constant $N$, depending
only on $K_{0}$, such that
for any $x,y\in\bR^{d}$, $\alpha_{\cdot}\in\frA$, 
$ \beta_{\cdot} 
\in\frB$, and $T\in(0,\infty)$,  
\begin{equation}
                                             \label{2.25.7}
E 
\sup_{t\leq T}|x_{t}^{\alpha_{\cdot}\beta_{\cdot}(x+y)} 
-x_{t}^{\alpha_{\cdot}\beta_{\cdot} x} |^{2}
\leq N|y|^{2} e^{NT}.
\end{equation}
Therefore
  for any $T\in(0,\infty)$ and $\kappa>0$
$$
\lim_{y\to 0}\sup_{x\in\bR^{d}}
\supsup_{\alpha_{\cdot}\in\frA\,\,\,\beta_{\cdot} 
\in\frB}
P 
(\sup_{t\leq T}|x_{t}^{\alpha_{\cdot}\beta_{\cdot}(x+y)} 
-x_{t}^{\alpha_{\cdot}\beta_{\cdot} x} |\geq\kappa)=0,
$$
which as in the case of \eqref{2.24.2} yields that
$$
\lim_{y\to 0}\sup_{x\in\bR^{d}}
|v(x+y)-v(x)|=0,
$$
that is $v$ is uniformly continuous in $\bR^{d}$.
 
Now, since the processes 
$x^{\alpha_{\cdot}\beta_{\cdot}x}_{t}(\varepsilon)$
are uniformly nondegenerate,
we know  (see the proof of Theorem 3.1 of \cite{Kr_13}) that
 \eqref{1.14.1} holds if we replace there
$v$, $x_{t}$, and $\phi_{t}$
 with $v^{\varepsilon}$, $x_{t}(\varepsilon)$,
 and $\phi_{t}(\varepsilon)$,
respectively. We want to pass to the limit
as $\varepsilon\downarrow0$ in the so modified 
\eqref{1.14.1}. By \eqref{2.24.2} the left-hand sides
will converge to $v(x)$.

It turns out that the limit of the right-hand sides
will not change if we replace back $v^{\varepsilon}$
with $v$. Indeed, the error of such replacement
is less than
$$
\sup_{\bR^{d}}
|v^{\varepsilon}-v|
\supsup_{\alpha_{\cdot}\in\frA\,\,\,\beta_{\cdot} 
\in\frB}E_{x}^{\alpha_{\cdot} \beta_{\cdot}}
\big[e^{-\psi_{\gamma}}+\int_{0}^{\gamma}\lambda_{t}
e^{-\psi_{t}}\,dt\big]=
\sup_{\bR^{d}}
|v^{\varepsilon}-v|.
$$
Hence, we reduced the proof of \eqref{1.14.1} to the proof
that the limit of
$$
\infsup_{\bbeta\in\bB\,\,\alpha_{\cdot}\in\frA}
E_{x}^{\alpha_{\cdot}\bbeta(\alpha_{\cdot})}\big[
v(x_{\gamma}(\varepsilon))e^{-\phi_{\gamma}(\varepsilon)-\psi_{\gamma}}
$$
\begin{equation}
                                                 \label{2.24.4}
+\int_{0}^{\gamma}
\{f( x_{t}(\varepsilon))+\lambda_{t}v(x_{t}(\varepsilon))\}
e^{-\phi_{t}(\varepsilon)-\psi_{t}}\,dt \big]
\end{equation}
equals the right-hand side of \eqref{1.14.1}.

As is easy to see the difference of \eqref{2.24.4}
and the right-hand side of \eqref{1.14.1} is less than
$I(\varepsilon)+J(\varepsilon)$, where
$$
I(\varepsilon)=\supsup_{\alpha_{\cdot}\in\frA\,\,\,\beta_{\cdot} 
\in\frB}E_{x}^{\alpha_{\cdot} \beta_{\cdot}}
\int_{0}^{\infty}\big|f( x_{t}(\varepsilon))
e^{-\phi_{t}(\varepsilon)}-
f( x_{t} )
e^{-\phi_{t} }\big|\,dt,
$$
$$
J(\varepsilon)=
\supsup_{\alpha_{\cdot}\in\frA\,\,\,\beta_{\cdot} 
\in\frB}E_{x}^{\alpha_{\cdot} \beta_{\cdot}}
\sup_{t\geq0}\big(|v(x_{t}(\varepsilon))
e^{-\phi_{t}(\varepsilon)}-v(x_{t})e^{-\phi_{t} }| 
\big) .
$$
Obviously, $I(\varepsilon)$ is less than the right-hand side 
of \eqref{2.24.5} and therefore tends to zero
as $\varepsilon\downarrow0$. The same is true
for $J(\varepsilon)$ which follows from the uniform
continuity of $v$ and $c$ and \eqref{2.24.1}.
The theorem is proved.

\begin{remark} 
 
It is unknown to the author
whether Theorem \ref{theorem 1.14.1} is still true
or not if we drop the assumption about the
existence of $\bar{w}_{t}$.

\end{remark}

\mysection{An auxiliary stochastic differential game
on a surface}
                                    \label{section 2.27.1}

Again the probability space here and the underlying filtration
of $\sigma$-fields are not necessarily the same
as in Section \ref{section 2.26.3} and in our applications
they indeed may be different. Therefore,
the following assumption is harmless for the purpose
of our applications.
\begin{assumption}
                               \label{assumption 2.25.3}
On $(\Omega,\cF,P)$ we are given four $d_{1}$-dimensional
and one $d+4$-dimensional
independent Wiener processes 
$ w^{1}_{t},...,w^{(4)}_{t},\bar{w}_{t}$, respectively,
which are Wiener processes relative to $\{\cF_{t}\}$.
\end{assumption}

We will work in
the space $\bR^{d}\times\bR^{4}=\{z=(x,y):x\in\bR^{d},
y\in\bR^{4}\}$.
 Set $\bar\Psi(x,y)=\Psi(x)-|y|^{2}$
and
in $\bR^{d}\times\bR^{4}$ consider the surface
$$
\Gamma=\{z:\bar\Psi(z)=0\}.
$$
The gradient of $\bar{\Psi}$ is not vanishing on $\Gamma$,
because the gradient of $\Psi$ is not vanishing on $\partial G$,
and, since $\bar\Psi\in C^{3}$, $\Gamma$ is a smooth surface
of class $C^{3}$. Obviously $\Gamma$ is closed and
bounded.

Denote by $D\Psi$ the gradient of $\Psi$ which we view as a column-vector
and set
$$
\hat
c^{\alpha\beta}(x)=-L^{\alpha\beta}\Psi(x)-c^{\alpha\beta}\Psi(x).
$$

Next, for $\alpha\in A,\beta\in B$, $z=(x,y)
\in \bR^{d}\times\bR^{4}$, and $i=1,...,4$ we define the functions 
$$
\bar{\sigma}^{\alpha\beta(i)}(z),\quad
 \bar{\sigma}^{\alpha\beta }(z),\quad
\bar{b}^{\alpha\beta(i)}(z),\quad
\bar{b}^{\alpha\beta }(z)
$$
in such a way that on $\Gamma$ they coincide with
$$
y^{i}\sigma^{\alpha \beta }(x ),\quad
(1/2)[D\Psi(x )]^{*}
\sigma^{\alpha \beta }(x ),\quad
-(1/2)y^{i} \hat 
c^{\alpha \beta }(x ),
$$
$$
|y|^{2}b^{\alpha \beta }
(x )+a^{\alpha \beta }
(x )D\Psi(x ),
$$
respectively, and are Lipschitz
continuous functions  of $z$ with compact support
 with Lipschitz constant and support
independent of $\alpha$ and $\beta$.

We also set
$$
\bar c^{\alpha\beta}(x,y)=-L^{\alpha\beta}\Psi(x)
$$
on $\Gamma$ and continue $\bar{c}^{\alpha\beta}(z)$ outside
$\Gamma$ in such a way that it is still Lipschitz
continuous in $z$ with Lipschitz constant independent
of $\alpha$ and $\beta$ and is greater than $ 1/2$ everywhere,
the latter being possible since $L^{\alpha\beta}\Psi\leq-1$
in $G$.

Next, we take $\alpha_{\cdot}\in\frA$, $\beta_{\cdot}
\in\frB$, $z=(x,y)\in \bR^{d}\times\bR^{4}$ and
define 
$$
z_{t}^{\alpha_{\cdot}\beta_{\cdot}z}
=(x,y)_{t}^{\alpha_{\cdot}\beta_{\cdot}z}
$$
 by means of the system
\begin{equation}
                                                   \label{2.25.1}
x_{t}=x+\int_{0}^{t} \bar\sigma^{\alpha_{s}\beta_{s}(i)}
(z_{s})\,dw^{(i)}_{s}+\int_{0}^{t}
\bar b^{\alpha_{s}\beta_{s}}
(z_{s}) \,ds,
\end{equation}
\begin{equation}
                                                   \label{2.25.2}
y^{i}_{t}=y^{i}+ \int_{0}^{t} 
\bar \sigma^{\alpha_{s}\beta_{s}}
(z_{s})\,dw^{(i)}_{s}+\int_{0}^{t}\bar 
b^{\alpha_{s}\beta_{s}}(z_{s})\,ds,
\end{equation}
$i=1,...,4$.

\begin{lemma}
                                             \label{lemma 2.25.1}
If $z\in \Gamma$, then $z_{t}^{\alpha_{\cdot}\beta_{\cdot}z}
\in \Gamma$ for all $t\geq0$ (a.s.) for any
$\alpha_{\cdot}\in\frA$ and $\beta_{\cdot}
\in\frB$ and $z_{t}^{\alpha_{\cdot}\beta_{\cdot}z}$
also satisfies the system
\begin{equation}
                                                   \label{2.11.1}
x_{t}=x+\int_{0}^{t}y^{i}_{s}\sigma^{\alpha_{s}\beta_{s}}
(x_{s})\,dw^{(i)}_{s}+\int_{0}^{t}
\big[|y_{s}|^{2}b^{\alpha_{s}\beta_{s}}
(x_{s})+2a^{\alpha_{s}\beta_{s}}
(x_{s})D\Psi(x_{s})\big]\,ds,
\end{equation}
\begin{equation}
                                                   \label{2.11.2}
y^{i}_{t}=y^{i}+(1/2)\int_{0}^{t}[D\Psi(x_{s})]^{*}
\sigma^{\alpha_{s}\beta_{s}}
(x_{s})\,dw^{(i)}_{s}-(1/2)\int_{0}^{t}y^{i}_{s}\hat 
c^{\alpha_{s}\beta_{s}}(x_{s})\,ds,
\end{equation}
$i=1,...,4$, in which one can replace $|y_{s}|^{2}$
with $\Psi(x_{s})$.
\end{lemma}

Proof. 
The system \eqref{2.11.1}-\eqref{2.11.2} has at least a local
solution before the solution explodes. However, the reader
will easily check by using It\^o's formula that
$d(\Psi(x_{t})-|y_{t}|^{2})=0$ and, since $\Psi$ 
is bounded from above,
$y_{t}$ cannot explode and  $x_{t}$ cannot explode either
since $\Psi(x)\to-\infty$ as $|x|\to\infty$.

In particular, if $(x,y)\in\Gamma$, then the solution
of \eqref{2.11.1}-\eqref{2.11.2} stays on $\Gamma$ for all
times. Then it satisfies \eqref{2.25.1}--\eqref{2.25.2},
and since the solution of the latter is unique,
the lemma is proved.

\begin{remark}
                                            \label{remark 2.27.1}
Observe that the process $z^{\alpha_{\cdot}\beta_{\cdot}z}_{t}$
is always a degenerate one and not only because
the coefficients of \eqref{2.25.1}-\eqref{2.25.2}
have compact support but also because, say, 
the diffusion   in \eqref{2.11.2} vanishes
when the $x$th component reaches (or just starts from)
 the maximum
point of $\Psi$, where $D\Psi=0$.
\end{remark}

Now we introduce a value function
$$
\bar v(z)=\infsup_{\bbeta\in\bB\,\,\alpha_{\cdot}\in\frA}
E^{\alpha_{\cdot}\bbeta(\alpha_{\cdot})}_{z}
\int_{0}^{\infty}f 
(x_{t})e^{-\bar\phi_{t}}\,dt,
$$
where
$$
\bar\phi_{t}^{\alpha_{\cdot}\beta_{\cdot}z}
=\int_{0}^{t}\bar c^{\alpha_{t}\beta_{t}}
(z_{s}^{\alpha_{\cdot}\beta_{\cdot}z})\,ds.
$$

Here is a fundamental fact relating the original differential
game in domain $G$, which is a domain
with boundary, with the one on $\Gamma$,
which is a closed manifold without boundary.

\begin{theorem}
                                         \label{theorem 2.11.1}
Suppose that $g\equiv0$. Then 
for $x\in G$ and $y\in\bR^{d}$ such that
$|y|^{2}=\Psi(x)$ we have $  \bar v(x,y)=v(x)/\Psi(x)$.

\end{theorem}

Proof. Fix $x\in G$ and $y\in\bR^{d}$ such that $|y|^{2}=\Psi(x)$
and take an $\varepsilon\in(0,\Psi(x))$.
Introduce, $z=(x,y)$ and
$$
\tau_{\varepsilon}^{\alpha_{\cdot}\beta_{\cdot}}=
\inf\{t>0:\Psi(x^{\alpha_{\cdot}\beta_{\cdot}z}_{t})
=\varepsilon\}.
$$
Then by Theorem \ref{theorem 1.14.1} (here we need
the existence of $\bar{w}_{t}$)
\begin{equation}
                                               \label{2.25.3}
\bar v(z)=
\infsup_{\bbeta\in\bB\,\,\alpha_{\cdot}\in\frA}
E_{z}^{\alpha_{\cdot}\bbeta(\alpha_{\cdot})}\big[
\bar v(z_{\tau_{\varepsilon}})e^{-\bar\phi_{\tau_{\varepsilon}}}
+\int_{0}^{\tau_{\varepsilon}}
 f ( x_{t}) e^{-\bar\phi_{t} }\,dt \big].
\end{equation}
By using It\^o's formula and Lemma \ref{lemma 2.25.1}
one easily sees that
$$
\Psi^{-1}(x^{\alpha_{\cdot}\beta_{\cdot}z}_{t})
\exp \big(-\int_{0}^{t}\hat c  ^{\alpha_{s}\beta_{s}}
 (x^{\alpha_{\cdot}\beta_{\cdot}z}_{s})\,ds\big)
$$
is a local martingale as long as it is well defined.
Since it is nonnegative it has bounded trajectories
implying that $\Psi (x^{\alpha_{\cdot}\beta_{\cdot}z}_{t})$
can never reach $0$ in finite time. Furthermore,
$$
E_{z}^{\alpha_{\cdot}\beta_{\cdot}}
e^{-\bar\phi_{\tau_{\varepsilon}}}=
E_{z}^{\alpha_{\cdot}\beta_{\cdot}}
e^{-\bar\phi_{\tau_{\varepsilon}}}I_{
\tau_{\varepsilon}<\infty}=\varepsilon
E_{z}^{\alpha_{\cdot}\beta_{\cdot}}\Psi^{-1}
(x_{\tau_{\varepsilon}})
e^{-\bar\phi_{\tau_{\varepsilon}}}I_{
\tau_{\varepsilon}<\infty}
$$
$$
\leq \varepsilon
E_{z}^{\alpha_{\cdot}\beta_{\cdot}}\Psi^{-1}
(x_{\tau_{\varepsilon}})
\exp \big(-\int_{0}^{\tau_{\varepsilon}}\hat c  ^{\alpha_{s}\beta_{s}}
 (x _{s})\,ds\big)I_{
\tau_{\varepsilon}<\infty}\leq \varepsilon \Psi^{-1}(x).
$$
This estimate is uniform with respect to $\alpha_{\cdot}$
and $\beta_{\cdot}$ and we conclude from \eqref{2.25.3}
that
\begin{equation}
                                               \label{2.25.4}
\bar v(z)=\lim_{\varepsilon\downarrow0}
\infsup_{\bbeta\in\bB\,\,\alpha_{\cdot}\in\frA}
E_{z}^{\alpha_{\cdot}\bbeta(\alpha_{\cdot})} \int_{0}^{\tau_{\varepsilon}}
 f ( x_{t}) e^{-\bar\phi_{t} }\,dt  .
\end{equation}

Next set
$$
\hat{w}_{t}^{\alpha_{\cdot}\beta_{\cdot}z}
=\int_{0}^{t} 
\Psi^{-1/2}(x_{s}^{\alpha_{\cdot}\beta_{\cdot}z})
(y_{s}^{\alpha_{\cdot}\beta_{\cdot}z})^{i}\,dw^{(i)}_{s},
$$
(recall that $\Psi (x_{s}^{\alpha_{\cdot}\beta_{\cdot}z})>0$
for all $s$).
 Since (a.s.) 
$$
|y_{s}^{\alpha_{\cdot}\beta_{\cdot}z}|^{2}
=\Psi(x_{s}^{\alpha_{\cdot}\beta_{\cdot}z})
$$
for all $s\geq0$, the process $
\hat{w}_{t}^{\alpha_{\cdot}\beta_{\cdot}z}$ is 
well defined and is a Wiener process.
Obviously it is control adapted in the terminology of \cite{Kr_14}.

Furthermore, 
$$
\int_{0}^{t}\Psi^{1/2}(x_{s}^{\alpha_{\cdot}\beta_{\cdot}z})
\sigma^{\alpha_{s}\beta_{s}}
(x_{s}^{\alpha_{\cdot}\beta_{\cdot}z})\,d
\hat{w}_{s}^{\alpha_{\cdot}\beta_{\cdot}z} 
$$
$$
=\int_{0}^{t} 
\sigma^{\alpha_{s}\beta_{s}}
(x_{s}^{\alpha_{\cdot}\beta_{\cdot}z})
(y_{s}^{\alpha_{\cdot}\beta_{\cdot}z})^{i}\,dw^{(i)}_{s}.
$$
We conclude that $x_{t}^{\alpha_{\cdot}\beta_{\cdot}z}$
satisfies the equation
$$
x_{t}=x+ \int_{0}^{t}\Psi^{1/2}(x_{s}) 
\sigma^{\alpha_{s}\beta_{s}}
(x_{s} )\,d
\hat{w}_{s}^{\alpha_{\cdot}\beta_{\cdot}z}
$$
\begin{equation}
                                            \label{2.12.1}
+\int_{0}^{t}
\big[\Psi(x_{s})b^{\alpha_{s}\beta_{s}}
(x_{s})+a^{\alpha_{s}\beta_{s}}
(x_{s})D\Psi(x_{s})\big]\,ds.
\end{equation}
Next,   define
$$
\quad r^{\alpha_{\cdot}\beta_{\cdot}}_{t}=
\Psi^{1/2}(x^{\alpha_{\cdot}\beta_{\cdot}z}_{t})I_{t\leq 
\tau_{\varepsilon}^{\alpha_{\cdot}\beta_{\cdot}}}
+I_{t>
\tau_{\varepsilon}^{\alpha_{\cdot}\beta_{\cdot}}}.
$$
Observe that $r^{\alpha_{\cdot}\beta_{\cdot}}_{t}$ is control adapted
($z$ is fixed) and for $t\leq \tau_{\varepsilon}^{\alpha_{\cdot}\beta_{\cdot}}$
the process $x^{\alpha_{\cdot}\beta_{\cdot}z}_{t}$ is a solution of
$$
x_{t}=x+ \int_{0}^{t}r^{\alpha_{\cdot}\beta_{\cdot}}_{s} 
\sigma^{\alpha_{s}\beta_{s}}
(x_{s} )\,d
\hat{w}_{s}^{\alpha_{\cdot}\beta_{\cdot}z}
$$
\begin{equation}
                                            \label{2.12.3}
+\int_{0}^{t}[r^{\alpha_{\cdot}\beta_{\cdot}}_{s}]^{2}
\big[  b^{\alpha_{s}\beta_{s}}
 +(\varepsilon\wedge\Psi^{-1})a^{\alpha_{s}\beta_{s}}
 D\Psi\big](x_{s})\,ds.
\end{equation}
Moreover,
for $t\leq \tau_{\varepsilon}^{\alpha_{\cdot}\beta_{\cdot}}$
$$
\bar\phi^{\alpha_{\cdot}\beta_{\cdot}z}_{t}
=\int_{0}^{t}[r^{\alpha_{\cdot}\beta_{\cdot}}_{s}]^{2}
(\varepsilon\wedge\Psi^{-1})\bar{c}
^{\alpha_{s}\beta_{s}}(x_{s}^{\alpha_{\cdot}\beta_{\cdot}z})\,ds.
$$
By Theorem 2.1 of \cite{Kr_14} (which, basically,
allows for random time changes and changes of probability measure
based on Girsanov's theorem)
$$
\bar v_{\varepsilon}(x):=\infsup_{\bbeta\in\bB\,\,\alpha_{\cdot}\in\frA}
E_{x}^{\alpha_{\cdot}\bbeta(\alpha_{\cdot})} 
\int_{0}^{\hat\tau_{\varepsilon}}
(\varepsilon\wedge\Psi^{-1})f ( \hat x_{t}) e^{-\hat\phi_{t} }\,dt, 
$$
where $\hat{x}_{t}^{\alpha_{\cdot}\beta_{\cdot}x}$
is a unique solution of
$$
x_{t}=x+ \int_{0}^{t}  
\sigma^{\alpha_{s}\beta_{s}}
(x_{s} )\,d
\hat{w}_{s}^{\alpha_{\cdot}\beta_{\cdot}z}
$$
$$
+\int_{0}^{t} 
\big[ b^{\alpha_{s}\beta_{s}}
 +(\varepsilon\wedge\Psi^{-1})a^{\alpha_{s}\beta_{s}}
 D\Psi\big](x_{s})\,ds,
$$
$$
\hat{\phi}_{t}^{\alpha_{\cdot}\beta_{\cdot}x}
=\int_{0}^{t} (\varepsilon\wedge\Psi^{-1})\bar{c}
^{\alpha_{s}\beta_{s}}
(\hat x_{s}^{\alpha_{\cdot}\beta_{\cdot}x})\,ds,
$$
$$
\hat{\tau}_{\varepsilon}^{\alpha_{\cdot}\beta_{\cdot}x}
=\inf\{t\geq 0:\Psi(
\hat{x}_{t}^{\alpha_{\cdot}\beta_{\cdot}x})\leq\varepsilon\}.
$$

Now it follows from \eqref{2.25.4} that
\begin{equation}
                                               \label{2.25.6}
\bar v(z)=\lim_{\varepsilon\downarrow0}
\bar v_{\varepsilon}(x).
\end{equation}

Also observe that by It\^o's formula, dropping
for simplicity of notation the indices
$\alpha_{\cdot},\beta_{\cdot},x$, we obtain
that for $t<\hat{\tau}_{\varepsilon}$
$$
\Psi^{-1}(\hat{x}_{t})e^{-\hat{\phi}_{t}}=
\Psi^{-1}(x)+\exp\big
[-\int_{0}^{t}\Psi^{-1}[D\Psi]^{*}
\sigma^{\alpha_{s}\beta_{s}}
(\hat x_{s} )\,d
\hat{w}_{s}
$$

$$
-\int_{0}^{t}\big[
\Psi^{-1}[D\Psi]^{*}b^{\alpha_{s}\beta_{s}}+
\Psi^{-2}[D\Psi]^{*}a^{\alpha_{s}\beta_{s}}D\Psi]
$$

$$
+
\Psi^{-1}\tr a^{\alpha_{s}\beta_{s}}D^{2}\Psi
-\Psi^{-1}L^{\alpha_{s}\beta_{s}}\Psi\big](\hat x_{s})\,ds\big].
$$
This result after obvious cancellations and 
introducing the notation
$$
\pi^{\alpha_{\cdot}\beta_{\cdot}x}_{t}=
(\varepsilon\wedge\Psi^{-1})[D\Psi]^{*}
\sigma^{\alpha_{t}\beta_{t}}
(\hat x^{\alpha_{\cdot}\beta_{\cdot}x}_{t} ) ,
$$
$$
\check\phi^{\alpha_{\cdot}\beta_{\cdot}x}_{t}
=\int_{0}^{t}c^{\alpha_{s}\beta_{s} } 
(\hat x^{\alpha_{\cdot}\beta_{\cdot}x}_{s})\,ds,
$$
$$
\psi^{\alpha_{\cdot}\beta_{\cdot}x}_{t}=
-\int_{0}^{t}\pi^{\alpha_{\cdot}\beta_{\cdot}x}_{s}\,
d\hat{w}_{s}^{\alpha_{\cdot}\beta_{\cdot}z}
-(1/2)\int_{0}^{t}|\pi^{\alpha_{\cdot}\beta_{\cdot}x}_{s}|^{2}
\,ds
$$
allows us to rewrite the definition of $\hat{v}_{\varepsilon}(x)$ as
$$
\bar v_{\varepsilon}(x)=\Psi^{-1}(x)
\infsup_{\bbeta\in\bB\,\,\alpha_{\cdot}\in\frA}
E_{x}^{\alpha_{\cdot}\bbeta(\alpha_{\cdot})} 
\int_{0}^{\hat\tau_{\varepsilon}}
  f ( \hat x_{t}) e^{-\check\phi_{t}-\psi_{t} }\,dt .
$$
Here by Theorem  2.1 of \cite{Kr_14} the right-hand side  is
equal to the expression
$$
\Psi^{-1}(x)\infsup_{\bbeta\in\bB\,\,\alpha_{\cdot}\in\frA}
E_{x}^{\alpha_{\cdot}\bbeta(\alpha_{\cdot})} 
\int_{0}^{ \tau_{\varepsilon}}
  f (  x_{t}) e^{- \phi_{t}  }\,dt 
$$
constructed on the probability space from
Section \ref{section 2.26.3} with $G_{\varepsilon}
=\{x:\Psi(x)>\varepsilon\}$ in place of $G$.
One shows that the limit as $\varepsilon\downarrow0$
of the last expression is $v(x)/\Psi(x)$
by repeating the proof of Theorem 2.2 of \cite{Kr_13_1}
given there in Section 6. After that
by coming back to \eqref{2.25.6}
one obtains the desired result. The theorem is proved.

This theorem allows us to make the first step
in proving approximation theorems by establishing
the Lipschitz continuity of $\bar v$ on $\Gamma$ away
from the equator.

 \begin{corollary}
                                        \label{corollary 2.25.1}
For any $\varepsilon>0$ there exists a constant
$N$ such that for any $z'=(x',y'),z''=(x'',y'')\in\Gamma$
satisfying $|y'|^{2},|y''|^{2}>\varepsilon$ we have
$$
|\bar v(z')-\bar v(z'')|\leq N|z'-z''|.
$$
\end{corollary}
  
Indeed, $\bar v(z')=v(x')/\Psi(x')$ and $
\bar v(z'')=v(x'')/\Psi(x'')$
and we know from \cite{Kr_14_1} (or from Remark 2.2
of \cite{Kr_14} and \cite{Tr_89}) that $v\in C^{0,1}_{loc}( G )$
(actually, $v$ belongs to a much better class).
Therefore, if $\Psi(x'),\Psi(x'')>\varepsilon$,
the difference $|\bar v(z')-\bar v(z'')|$ is less than
a constant times $|x'-x''|\leq|z'-z''|$.

To establish the  
Lipschitz continuity of $\bar v$ on the whole of $\Gamma$
we need the following.
\begin{lemma}
                                             \label{lemma 2.26.1}
(i) There is a constant $N_{0}$,
depending only on the Lipschitz constants
of the coefficients of \eqref{2.25.1}-\eqref{2.25.2},
such that for any $z',z''\in \bR^{d}\times\bR^{4}$,
$\alpha_{\cdot}\in\frA$, and $\beta_{\cdot}
\in\frB$
the process
$$
|z_{t}^{\alpha_{\cdot}\beta_{\cdot}z'}-
z_{t}^{\alpha_{\cdot}\beta_{\cdot}z''}|^{2}e^{-2N_{0}t}
+\int_{0}^{t}|z_{t}^{\alpha_{\cdot}\beta_{\cdot}z'}-
z_{s}^{\alpha_{\cdot}\beta_{\cdot}z''}|^{2}
e^{-2N_{0}s}\,ds
$$
is a supermartingale.

(ii) There exists an $\varepsilon>0$ such that
if $z=(x,y)\in\Gamma$ and $|y|^{2}\leq \varepsilon$, then
for any $\alpha_{\cdot}\in\frA$ and $\beta_{\cdot}
\in\frB$
$$
E^{\alpha_{\cdot}\beta_{\cdot}}_{z}e^{2N_{0}\tau_{2\varepsilon}}
\leq \frac{1}{\cos1}.
$$
\end{lemma}

Proof. Assertion (i) is easily obtained after computing
the stochastic differential of the process in question.

To prove (ii), observe that 
$|D\Psi|\geq1$ on $\partial G$ and hence for a sufficiently small
$\varepsilon>0$ we have  $|D\Psi|\geq1/2$ if 
$\Psi\in[0,2\varepsilon]$.
  In that case also
$$
\nu^{\alpha\beta}:=
a^{\alpha\beta}_{ij}(D_{i}\Psi)D_{j}\Psi\geq \delta/4.
$$
Next, denote $\lambda=(2\varepsilon)^{-1/2}$ and note that by
  It\^o's formula,
 dropping the indices $\alpha_{\cdot}$, $\beta_{\cdot}$, and $z$,
one obtains
$$
d\big[e^{2N_{0}t}\cos\lambda|y_{t}|\big]=e^{2N_{0}t}
(\lambda|y_{t}|/2)\hat{c}^{\alpha_{t}\beta_{t}}
(x_{t})\sin\lambda|y_{t}|\,dt
$$
$$
-e^{2N_{0}t}\bigg[\frac{3}{4}\frac{\lambda\sin\lambda|y_{t}|}{|y_{t}|}
\nu^{\alpha_{t}\beta_{t}}(x_{t})
+\frac{\lambda^{2}}{4}\nu^{\alpha_{t}\beta_{t}}(x_{t})\cos\lambda|y_{t}|-
2N_{0}\cos\lambda|y_{t}|\bigg]\,dt+dm_{t},
$$
where $m_{t}$ is a martingale starting from zero.
For $t\leq \tau_{2\varepsilon}$ the first term on
the right
is dominated by  $N_{1}e^{2N_{0}t}\,dt$, where $N_{1}$
is a constant,
since $\hat{c}^{\alpha\beta}(x)$ is bounded.
It is seen that reducing $\varepsilon$ if necessary
so that $\lambda=(2\varepsilon)^{-1/2}$ satisfies
$$
\frac{\lambda^{2}}{16}\delta\cos 1-
2N_{0}\cos 1\geq N_{1} ,
$$
we have for $t\leq \tau_{2\varepsilon}$ that
$$
d\big[e^{2N_{0}t}\cos\lambda|y_{t}| \big]\leq dm_{t}.
$$
It follows that
$$
\cos 1
E^{\alpha_{\cdot}\beta_{\cdot}}_{z}
 e^{2N_{0}\tau_{2\varepsilon}}
\leq E^{\alpha_{\cdot}\beta_{\cdot}}_{z}
\big[e^{2N_{0}\tau_{2\varepsilon}}
\cos\lambda|y_{\tau_{2\varepsilon}}|\big]
 \leq 1,
$$
and the lemma is proved.

\begin{theorem}
                                           \label{theorem 2.25.1}
There exists a constant
$N$ such that for any $z'=(x',y'),z''=(x'',y'')\in\Gamma$
 we have
\begin{equation}
                                                 \label{2.27.1}
|\bar v(z')-\bar v(z'')|\leq N|z'-z''|.
\end{equation}
\end{theorem}

Proof. Take  $\varepsilon>0$ from Lemma \ref{lemma 2.26.1}
 and fix $z'=(x',y'),z''=(x'',y'')\in\Gamma$ such that
$\Psi(x'),\Psi(x'')\leq 2\varepsilon$. Then
on the basis of Theorem \ref{theorem 1.14.1} write
$$
\bar v(z)=\infsup_{\bbeta\in\bB\,\,\alpha_{\cdot}\in\frA}
E_{z}^{\alpha_{\cdot}\bbeta(\alpha_{\cdot})}\big[
\bar v(z_{\gamma})e^{-\bar\phi_{\gamma} }
+\int_{0}^{\gamma}
 f^{\alpha_{t}\beta_{t}}( x_{t}) e^{-\bar\phi_{t} }\,dt \big],
$$
where
$$
\gamma^{\alpha_{\cdot}\beta_{\cdot}z}
=\tau_{2\varepsilon}^{\alpha_{\cdot}\beta_{\cdot}z'}
\wedge \tau_{2\varepsilon}^{\alpha_{\cdot}\beta_{\cdot}z''}.
$$

Next, fix $\alpha_{\cdot}\in\frA$ and $\beta_{\cdot}\in\frB$
and denote
$$
\tau'=\tau_{2\varepsilon}^{\alpha_{\cdot}\beta_{\cdot}z'},
\quad\tau''=\tau_{2\varepsilon}^{\alpha_{\cdot}\beta_{\cdot}z''},
\quad \gamma=\tau'\wedge\tau'', 
$$
$$
z'_{t}=z^{\alpha_{\cdot}\beta_{\cdot}z'}_{t},\quad
z''_{t}=z^{\alpha_{\cdot}\beta_{\cdot}z''}_{t},\quad
\bar{\phi}_{t}'=\bar{\phi}^{\alpha_{\cdot}\beta_{\cdot}z'}_{t},\quad
\bar{\phi}_{t}''=\bar{\phi}^{\alpha_{\cdot}\beta_{\cdot}z''}_{t}.
$$
Observe that 
$$
E\big|\bar v(z'_{\gamma  })
e^{-\bar\phi'_{\gamma }}-\bar v (z''_{\gamma })e^{-\bar\phi''
_{\gamma }}\big|\leq I_{1}+I_{2},
$$
where
$$
I_{1}=E |\bar v(z'_{\gamma  })
 -\bar v (z''_{\gamma })|,\quad I_{2}=
NE\int_{0}^{\gamma}|\bar{c}^{\alpha_{t}\beta_{t}}(x'_{t})
-\bar{c}^{\alpha_{t}\beta_{t}}(x''_{t})|\,dt.
$$

Below by $N$ we denote various constants independent of
$z',z'',\alpha_{\cdot}$, and $\beta_{\cdot}$.
By Corollary \ref{corollary 2.25.1} and Lemma
\ref{lemma 2.26.1}
$$
E |\bar v(z'_{\gamma  })
 -\bar v (z''_{\gamma })|I_{\Psi(x'_{\gamma}),\Psi(x''_{\gamma})
\geq\varepsilon}
\leq NE|z'_{\gamma  }-z''_{\gamma  }|
$$
$$
\leq NE^{1/2}|z'_{\gamma  }-z''_{\gamma  }|^{2}e^{-2N_{0}\gamma}
E^{1/2}e^{ 2N_{0}\gamma}\leq N|z'-z''|.
$$
 Furthermore,
$$
E |\bar v(z'_{\gamma  })
 -\bar v (z''_{\gamma })|I_{\Psi(x'_{\gamma})<\varepsilon,
\Psi(x''_{\gamma})
\geq\varepsilon}\leq N
E  I_{\Psi(x'_{\gamma})<\varepsilon,
\Psi(x''_{\gamma})
=2\varepsilon}
$$
$$
\leq\varepsilon^{-1}E|\Psi(x'_{\gamma})-\Psi(x''_{\gamma})|
\leq N|z'-z''|.
$$
Similarly,
$$
E |v(z'_{\gamma  })
 -v (z''_{\gamma })|I_{\Psi(x'_{\gamma})\geq\varepsilon,
\Psi(x''_{\gamma})
<\varepsilon}\leq N|z'-z''|
$$
and we conclude that $I_{1}\leq N|z'-z''|$.

Also by Lemma
\ref{lemma 2.26.1}
$$
I_{2}\leq NE\int_{0}^{\gamma}|x'_{t}-x''_{t}|\,dt
$$
$$
\leq NE^{1/2}\int_{0}^{\gamma}|x'_{t}-x''_{t}|^{2}e^{-2N_{0}
t}\,dtE^{1/2}\int_{0}^{\gamma}e^{ 2N_{0}
t}\,dt
$$
\begin{equation}
                                                 \label{2.27.3}
\leq N|z'-z''|E^{1/2} e^{ 2N_{0}\gamma}
\leq  N|z'-z''|.
\end{equation}
Hence,
\begin{equation}
                                                 \label{2.27.2}
E\big|\bar v(z'_{\gamma  })
e^{-\bar\phi'_{\gamma }}-\bar v (z''_{\gamma })e^{-\bar\phi''
_{\gamma }}\big| \leq  N|z'-z''|.
\end{equation}

Next, by using the inequalities $|e^{-a}-e^{-b}|\leq e^{-t}|a-b|$
valid for $a,b\geq t$ and $|ab-cd|\leq |b|\cdot|a-c|+|c|\cdot
|b-d|$ we obtain
$$
 \int_{0}^{\gamma}\big|
 f^{\alpha_{t}\beta_{t}}( x'_{t}) e^{-\bar\phi'_{t} }
- f^{\alpha_{t}\beta_{t}}( x''_{t}) e^{-\bar\phi''_{t} }\big|\,dt  
$$
$$
\leq\int_{0}^{\gamma}\big[|
f^{\alpha_{t}\beta_{t}}( x'_{t})-
f^{\alpha_{t}\beta_{t}}( x''_{t})|+
e^{-t}\int_{0}^{t}|
c^{\alpha_{s}\beta_{s}}( x'_{s})-
c^{\alpha_{s}\beta_{s}}( x''_{s})|\,ds\big]\,dt
$$
$$
\leq\int_{0}^{\gamma}\big[|
f^{\alpha_{t}\beta_{t}}( x'_{t})-
f^{\alpha_{t}\beta_{t}}( x''_{t})|+
 |
c^{\alpha_{t}\beta_{t}}( x'_{t})-
c^{\alpha_{t}\beta_{t}}( x''_{t})| \big]\,dt
$$
$$
\leq N\int_{0}^{\gamma}|x'_{t}-x''_{t}|\,dt.
$$
This along with \eqref{2.27.2} and \eqref{2.27.3} shows that
\eqref{2.27.1} holds if $\Psi(x'),\Psi(x'')\leq2\varepsilon$.

If $\Psi(x')\geq2\varepsilon$
and $\Psi(x'')\leq\varepsilon$, then $\varepsilon\leq
\Psi(x')-\Psi(x'')\leq N|x'-x''|$ and then certainly
\eqref{2.27.1} holds. The same happens if
$\Psi(x'')\geq2\varepsilon$
and $\Psi(x')\leq\varepsilon$

The remaining cases where 
$\Psi(x')\geq2\varepsilon$
and $\Psi(x'')\geq\varepsilon$ or
$\Psi(x'')\geq2\varepsilon$
and $\Psi(x')\geq\varepsilon$
are taken care
of by Corollary \ref{corollary 2.25.1}.
The theorem is proved. 

\mysection{Proof of Theorem \protect\ref{theorem 1.14.01}}
                                     \label{section 4.9.2}

Denote $A_{1}=A$ and let
let $A_{2}$ be a 
separable metric space having no common points with $A_{1}$.
Assume that on $A_{2}\times B\times\bR^{d}$
we are given bounded 
continuous
functions $\sigma^{\alpha}=\sigma^{\alpha\beta}$,
$b^{\alpha}=b^{\alpha\beta}$, $c^{\alpha}=c^{\alpha\beta}$  
(independent of  $x$ and $\beta$), 
and $f^{\alpha\beta}\equiv0$ 
satisfying  the assumptions in Section \ref{section 2.26.3}
  perhaps with   different constants $\delta $ 
and $K_{0}$. Actually, the concrete values
of these constants never played any role, so that
we can take them to be the same here and in
Section \ref{section 2.26.3} (take the largest
$K_{0}$ as a new $K_{0}$ and the smallest...).
We made this comment to be able to use the same 
function $\Psi$ here as in Section \ref{section 2.26.3}.
 
Define
$$
\hat{A}=A_{1}\cup A_{2}.  
$$

Then we introduce $\hat{\frA}$ as the set of progressively measurable
$\hat{A}$-valued processes and $\hat{\bB}$ as the set
of $\frB$-valued functions $ \bbeta(\alpha_{\cdot})$
on $\hat{\frA}$ such that,
for any $T\in[0,\infty)$ and any $\alpha^{1}_{\cdot},
\alpha^{2}_{\cdot}\in\hat{\frA}$ satisfying
$$
P(  \alpha^{1}_{t}=\alpha^{2}_{t} 
 \quad\text{for almost all}\quad t\leq T)=1,
$$
we have
$$
P(  
 \bbeta_{t}(\alpha^{1}_{\cdot})=\bbeta_{t}(\alpha^{2}_{\cdot}) 
 \quad\text{for almost all}\quad t\leq T)=1.
$$

Next, take a constant $K\geq0$ and set 
$$
v_{K}(x)=\infsup_{\bbeta\in\hat{\bB}\,\,\alpha_{\cdot}\in\hat{\frA}}
v^{\alpha_{\cdot}\bbeta(\alpha_{\cdot})}_{K}(x),
$$
where
$$
v^{\alpha_{\cdot}\beta _{\cdot} }_{K}(x)=
E_{x}^{\alpha_{\cdot}\beta _{\cdot} }
\big[ \int_{0}^{\tau}
 f _{K} ( x_{t})e^{-\phi_{t} }\,dt 
+g(x_{\tau})e^{-\phi_{\tau} }\big]
$$
$$
 f^{\alpha\beta}_{K}( x)=f^{\alpha\beta}( x)-KI_{\alpha\in A_{2}}.
$$

As is explained in Section 6 of \cite{Kr_14} 
there is a set $A_{2}$ and other objects mentioned
above such that $u_{K}=v_{K}$ in $G$. Observe that
$|v-v_{K}|=|(v-g)-(v_{k}-g)|$ and
since $g\in C^{3}$ we can transform $v-g$ and 
$v_{k}-g$ 
by using It\^o's formula. Then we see that
$$
v(x)-g(x)=
\infsup_{\bbeta\in\bB\,\,\alpha_{\cdot}\in\frA}
E_{x}^{\alpha_{\cdot}\bbeta(\alpha_{\cdot})}
 \int_{0}^{\tau}
[Lg+ f] ( x_{t})e^{-\phi_{t}
 }\,dt ,
$$
where 
\begin{equation}
                                                    \label{2.27.06}
L^{\alpha\beta}g(x)+f^{\alpha\beta}(x),
\end{equation}
 $\alpha\in A$,
$\beta\in B$, $x\in\bR^{d}$, now plays
the role of a new $f^{\alpha\beta}(x)$ and possesses the same 
regularity properties as the old one. Also
$$
v_{K}(x)-g(x)=
\infsup_{\bbeta\in\hat{\bB}\,\,\alpha_{\cdot}\in\hat{\frA}}
E_{x}^{\alpha_{\cdot}\beta _{\cdot} }
 \int_{0}^{\tau}
[Lg+ f _{K} ]( x_{t})e^{-\phi_{t} }\,dt . 
$$
We see that, by replacing the original $f^{\alpha\beta}(x)$
with expression \eqref{2.27.06} (for
$\alpha\in \hat A$,
$\beta\in B$, $x\in\bR^{d}$) we reduce the proof of the theorem
to the proof that
\begin{equation}
                                                    \label{2.27.6}
|v-v_{K}|\leq N\Psi/K
\end{equation}
in $G$ for $K\geq1$ if $g\equiv0$. The only additional
change with regard to the setting in the beginning
of the section is that the new $f^{\alpha\beta}(x)$
generally is not zero when $\alpha\in A_{2}$.
With this in mind we proceed further assuming that
$$
g\equiv0.
$$

Now, if necessary, we pass to a
different complete probability space 
$(\bar \Omega,\bar P,\bar \cF)$
with an increasing filtration 
$\{\bar \cF_{t},t\geq0\}$  
of $\sigma$-fields $\bar \cF_{t}\subset \bar \cF $ such that
each $\bar \cF_{t}$ is complete with respect to $\bar \cF,\bar P$.
We can find such a space so that it  carries
  four $d_{1}$-dimensional
and one $d+4$-dimensional
independent Wiener processes 
$ w^{1}_{t},...,w^{(4)}_{t},\bar{w}_{t}$,
which are Wiener processes relative to 
$\{\bar \cF_{t},t\geq0\}$.
After that we repeat the constructions in 
Section \ref{section 2.27.1}
replacing there $A$ with $\hat A$
(now, of course, $\alpha_{t}$ and $\beta_{t}$
are $\hat A$- and $B$-valued functions, respectively,
  defined on $\bar{\Omega}$). Fix an element $\alpha^{*}
\in A_{1}$ and define a projection operator
$p:\hat A\to A_{1}$ by $p\alpha=\alpha$ if $\alpha\in A_{1}$
and $p\alpha=\alpha^{*}$ if $\alpha\in A_{2}$

Next, we introduce  value functions
$$
\bar v(z)=\infsup_{\bbeta\in\bB\,\,\alpha_{\cdot}\in\frA}
\bar E^{p\alpha_{\cdot}\bbeta(p\alpha_{\cdot})}_{z}
\int_{0}^{\infty}f 
(x_{t})e^{-\bar\phi_{t}}\,dt,
$$
$$
\bar v_{K}(z)=\infsup_{\bbeta\in\bB\,\,\alpha_{\cdot}\in\frA}
\bar E^{\alpha_{\cdot}\bbeta(\alpha_{\cdot})}_{z}
\int_{0}^{\infty}f_{K} 
(x_{t})e^{-\bar\phi_{t}}\,dt.
$$
We keep the notation $\bar v(z)$ the same as in 
Section \ref{section 2.27.1}
since these two objects coincide
if the probability space, filtration, and the Wiener
processes coincide, because the range of $p\alpha$
is just $A$. They coincide even if 
the probability space, filtration, and the Wiener
processes are different owing to Theorem 2.1 of
\cite{Kr_14}.

Observe that obviously $\bar v_{K}\geq \bar v$ and now
in light of Theorem \ref{theorem 2.11.1} to prove 
\eqref{2.27.6}  it suffices to prove that
on $\Gamma$
\begin{equation}
                                                  \label{2.27.4}
\bar{v}_{K}\leq \bar v+N/K
\end{equation}
for $K\geq 1$ with $N$ being a constant.

We are, basically, going to repeat the proof
of Theorem 2.4 of \cite{Kr_14_1} given there
in Section  10  for the uniformly nondegenerate case.
In this connection see Remark \ref{remark 2.27.1}.

Define
$$
d_{K}=\sup_{\Gamma}(\bar v_{K}-\bar v),\quad\lambda=\supsup_{\alpha\in
\hat A\,\,\,\beta\in B}\sup_{ z\in\bR^{d+4}}
\bar c^{\alpha\beta}(z)
$$
and denote by $z$ a point in $\Gamma$ at which $d_{K}$ is attained.

By the dynamic programming principle (Theorem \ref{theorem 1.14.1})
$$
\bar v_{K}(z)=\infsup_{\bbeta\in \bB \,\,\alpha_{\cdot}\in\ \frA }
\bar E_{z}^{\alpha_{\cdot}\bbeta(\alpha_{\cdot})}
\big[\bar v_{K}(z_{1})e^{ -\lambda }
+\int_{0}^{1}\{f _{K}
+(\lambda-\bar c  ) \bar v_{K} \}(z_{t})
e^{ -\lambda t}\,dt\big].
$$
Observe that
$$
e^{-\lambda}+\int_{0}^{1}[\lambda-\bar c^{\alpha_{t}\beta_{t}}
(z^{\alpha_{\cdot}\beta_{\cdot}z}_{t})]e^{-\lambda t}\,dt
\leq e^{-\lambda}+\int_{0}^{1}(\lambda-1/2)
e^{-\lambda t}\,dt
=:\kappa<1.
$$
Hence,
$$
\bar v_{K}(z)\leq
\infsup_{\bbeta\in \bB \,\,\alpha_{\cdot}\in\ \frA }
\bar E_{z}^{\alpha_{\cdot}\bbeta(\alpha_{\cdot})}
\big[\bar v (z_{1})e^{ -\lambda }
+\int_{0}^{1}\{f _{K}+(\lambda-\bar c  ) 
\bar v  \}(z_{t})
e^{ -\lambda t}\,dt\big]+\kappa d_{K}.
$$

Now take a sequence    $\bbeta^{n}\in \bB$ such that
\begin{equation}
                                                  \label{6.3.3}
\bar v(z)\geq \sup_{\alpha_{\cdot}\in\frA}
\bar E^{p\alpha_{\cdot}\bbeta^{n}(p\alpha_{\cdot})}_{z}
\big[\int_{0}^{1}(f 
+(\lambda-\bar c)\bar v)(z_{t})e^{-\lambda t}\,dt
+e^{-\lambda}\bar v(z_{1})\big]-1/n.
\end{equation}
Then find $\alpha_{\cdot}^{n}\in\frA$
such that
$$
\bar v_{K}(z)\leq 
\bar E_{z}^{\alpha^{n}_{\cdot}\bbeta^{n}( p\alpha^{n}_{\cdot})}
\big[v (z_{1})e^{ -\lambda }
+\int_{0}^{1}\{f _{K}+(\lambda-\bar c  )
\bar  v  \}(z_{t})
e^{ -\lambda t}\,dt\big]+\kappa d_{K}+1/n
$$
\begin{equation}
                                                  \label{6.3.4}
=\bar E_{z}^{\alpha^{n}_{\cdot}\bbeta^{n}(p\alpha^{n}_{\cdot})}
\big[v (z_{1})e^{ -\lambda }
+\int_{0}^{1}\{f  +(\lambda-\bar c  ) \bar v  \}(z_{t})
e^{ -\lambda t}\,dt\big]
\end{equation}
$$
-KR_{n} 
+\kappa d_{K}+1/n,
$$
where
$$
R_{n}=\bar E\int_{0}^{1}e^{-\lambda t}I_{\alpha^{n}_{t}\in A_{2}}\,dt.
$$

By Lemma 5.3 of \cite{Kr_13} for any $\alpha_{\cdot}
\in \frA$ and $\beta_{\cdot}\in\frB$  we have
$$
\bar E\sup_{t\leq 1}|z_{t}^{p\alpha_{\cdot}\beta_{\cdot}z}
-z_{t}^{ \alpha_{\cdot}\beta_{\cdot}z}|\leq
N\big(\bar E^{\alpha_{\cdot}\beta_{\cdot}}_{z}\int_{0}^{1}
 I_{\alpha^{n}_{t}\in A_{2}}\,dt\big)^{1/2},
$$
where the constant $N$ depends only on $K_{0}$ 
  and $d$.
We use this and
since $\bar c,f,\bar v$ are Lipschitz continuous on $\Gamma$,
 we get from 
\eqref{6.3.4} and \eqref{6.3.3}
$$
\bar v_{K}(z)+(K-N_{0})R_{n}\leq 
 E_{z}^{p\alpha^{n}_{\cdot}\bbeta^{n}(p\alpha^{n}_{\cdot})}
\big[v (z_{1})e^{ -\lambda }
+\int_{0}^{1}\{f  +(\lambda-\bar c  ) \bar v  \}(z_{t})
e^{ -\lambda t}\,dt\big]
$$
$$ 
+\kappa d_{K}+1/n+NR^{1/2}_{n}\leq
\bar  v(z)+\kappa d_{K}+2/n+NR^{1/2}_{n},
$$
where the constant $N_{0}$ depends only on the supremums of $\bar c$,
$|\bar v|$, and $|f|$.
Hence
\begin{equation}
                                                  \label{6.3.6}
\bar v_{K}(z)-\bar v(z)-\kappa d_{K}+(K-N_{0})R_{n}
\leq 2/n+NR^{1/2}_{n}.
\end{equation}
Here $\bar v_{K}(z)-\bar v(z)-\kappa d_{K}=(1-\kappa)d_{K}$
which is nonnegative.
It follows that  
$$
(K-N_{0})R_{n}\le 2/n+NR^{1/2}_{n},
$$
which for $K\geq 2N_{0}+1$ implies that $K
R_{n}\le 4/n+NR^{1/2}_{n}$, so that, if $KR_{n}\geq
8/n$, then $KR_{n}\leq NR_{n}^{1/2}$ and
$R_{n}\leq N/K^{2}$. Thus,
$$
R_{n}\leq 8/(nK)+N/K^{2},
$$
which after coming back to \eqref{6.3.6} finally yields
$$
(1-\kappa)d_{K} 
\leq 2/n+N/\sqrt{n}+N/K.
$$
After letting $n\to\infty$ we obtain \eqref{2.27.4}
for $K\geq 2N_{0}+1$. For smaller $K$ the estimate
holds just because $\bar v$ and $\bar v_{K}$ are bounded.
The theorem is proved.

\end{document}